\documentclass[10pt]{amsart}
\usepackage{a4}
\usepackage{amssymb}
\usepackage{amscd}
\usepackage{array}
\usepackage{booktabs}
\usepackage{multirow}
\usepackage{hhline}
\usepackage[all]{xy}

\newtheorem{thm}{Theorem}

\newtheorem{cor}[thm]{Corollary}
\newtheorem{lem}[thm]{Lemma}
\newtheorem{prop}[thm]{Proposition}
\newtheorem{defn}[thm]{Definition}

\newtheorem{rem}[thm]{Remark}

\numberwithin{thm}{section}
\numberwithin{equation}{section}
\newcommand{\norm}[1]{\left\Vert#1\right\Vert}
\newcommand{\abs}[1]{\left\vert#1\right\vert}

\newcommand{\la}{\langle}
\newcommand{\ra}{\rangle}

\newcommand{\Hi}{\mathcal{H}}

\newcommand{\K}{\mathcal{K}}

\newcommand{\M}{\mathcal{M}}
\newcommand{\m}{\textbf{mod}}
\newcommand{\n}{\mathbb{N}}

\newcommand{\prt}{\widehat{\otimes}}

\begin{document}
\title{Projectivity of modules over Segal algebras}
\author{Brian E. Forrest, Hun Hee Lee, Ebrahim Samei}

\address{Brian E. Forrest: Department of Pure Mathematics, Faculty of Mathematics, University of Waterloo,
200 University Avenue West, Waterloo, Ontario, Canada N2L 3G1}
\address{Hun Hee Lee: Department of Mathematics, Chungbuk National University, 410 Sungbong-Ro, Heungduk-Gu, Cheongju 361-763, Korea}
\address{Ebrahim Samei: Department of Mathematics and Statistics, University of Saskatchewan, 106 Wiggins Road Saskatoon, SK S7N 5E6 CANADA}


\begin{abstract}
In this paper we will study the projetivity of various natural modules associated to operator Segal algebras 
of the Fourier algebra of a locally compact group. In particular, we will focus on the question of 
identifying when such modules will be projective in the category of operator spaces.
Projectivity often implies that the underlying group is discrete or even finite. 
We will also look at the projectivity for modules of $A_{cb}(G)$, the closure of $A(G)$ in the space of its completely bounded 
mutipliers. Here we give an evidence to show that weak amenability of $G$ plays an important role. 

\end{abstract}

\maketitle

\section{Introduction}
Recently, H. G. Dales and M. E. Polyakov (\cite{DP04}) gave a detailed study of the homological 
properties of modules over the group algebra of a locally compact group. In their work
they focused primarily on the question of whether or not certain natural left $L^1(G)$-modules are respectively 
\textit{projective}, \textit{injective}, or \textit{flat}. They were able to show, for example, that 
when viewed as the dual left-module of $L^{1}(G)$, $L^{\infty}(G)$ is projective precisely when the 
group $G$ is finite. In stark contrast, they proved that $L^{\infty}(G)$  is always injective. They also showed that the measure algebra $M(G)$
is projective precisely when $G$ is discrete, while in this case injectivity is equivalent to the group $G$ being 
amenable. 

It is well known that for abelian locally compact groups that classical Fourier transform 
identifies $L^{1}(G)$ with a commutative Banach algebra 
$A(\widehat {G})$ of continuous functions on the dual group $\widehat G$ called the Fourier algebra. If $G$ is not abelian, 
then the classical method for defining the Fourier algebra is no longer available. However, 
in \cite{Eymard} P. Eymard succeeded in defining the Fourier algebra $A(G)$ for any locally compact group $G$. 
Just as was the case for an abelian group, $A(G)$ is a commutative Banach algebras of continuous functions on $G$. 

Motivated by the earlier work of Dales and Polyakov, in the prequel to this paper \cite{FLS}, we studied the projectivity of various natural 
$A(G)$-modules when $G$ is a locally compact, not necessarily abelian group. We found that to do so it was necessary to 
incorporate the Fourier algebra's natural operator space structure into the investigation. In doing so we were able to 
obtain operator space analogues of most of Dales and Polyakov's results concerning the projectivity of modules over $L^{1}(G)$, as well as to 
deduce related results for additional modules. 

In this paper, we continue from where we left of in \cite{FLS} to further our study of projective modules arising from the 
Fourier algebra. However instead of focusing on $A(G)$ as our base algebra, we will look first at 
certain operator Segal algebras of $A(G)$, namely the Lebesgue-Fourier algebra, $S^1A(G):=L^1(G)\cap A(G)$ and the 
Feichtinger's algebra $S_0(G)$. Both of these are completely contractive Banach algebras, that are dense 
two-sided ideals in $A(G)$ with natural operator space structures that make them essential completely
contractive $A(G)$-modules. Both of these algebras have received significant attention lately. They each have 
a rich structure of their own but can also be used to reveal new insight into the nature of $A(G)$ itself, as well 
as that of the underlying group $G$. (See \cite{FSW} and \cite{Spr2} for a more detailed account of theses algebras.)

As we did in \cite{FLS}, we will show that discreteness plays an imprortant role in identifying projective modules 
for both $S^1A(G)$ and $S_0(G)$. This time, we also see that compactness plays an equal role. As such we find that 
many natural spaces are projective when viewed as modules over $S^1A(G)$ or $S_0(G)$ if and only if $G$ is finite. There are 
exceptions to this rule noteably $A(G)$, which is projective when $G$ is compact, and $L^1(G)$, which is projective if and only if $G$ is discrete.
This latter result reflects the natural duality between $A(G)$ and $L^1(G)$, though at this point it remains open 
whether or not the projectivity of $A(G)$ as either an $S^1A(G)$ or an $S_0(G)$-module forces $G$ to be compact. 

In the latter part of the paper, we will change directions once more. This time instead of looking at an operator 
Segal algebra of $A(G)$, we will instead focus our attention on another algebra $A_{cb}(G)$ in which $A(G)$ is 
itself the operator Segal algebra. 

In general, $A_{cb}(G)$ is the closure of $A(G)$ in the space of its completeley bound multipliers. However, for 
amenable groups it is exaclty $A(G)$. In contrast, when $G$ is non-amenable,  $A(G)$ is a proper operator 
Segal algebra of $A_{cb}(G)$.
In the case, that $G$ is weakly amenable, that is $A(G)$ has an approximate identity that is bounded in the cb-multiplier norm, 
$A_{cb}(G)$ is known to have many of the properties of the Fourier algebra of an amenable group. For example, when $G$ 
is the free group $\mathbb{F}_2$ on two generators, $G$ is weakly amenable and $A_{cb}(G)$ is known to be operator 
amenable (\cite{FRS07}). In contrast, $A(G)$ is operator amenable if and only if $G$ is amenable. 

We will show that as was the case, with the Fourier algebra, for many natural modules of $A_{cb}(G)$ to be projective 
we again require the underlying group to be discrete. While for some modules we even require the group to be finite, 
for groups like $\mathbb{F}_2$, we are able to establish projectivity for a number of important modules.

\section{Preliminaries}

\subsection{Harmonic analysis} Let $G$ be a locally compact group with a fixed left Haar measure $dx$.
Let $L^{1}(G)$ denote the group algebra of $G$. Then $L^{1}(G)$ is an involutive Banach algebra under convolution.

Let $\Sigma _{G}$ be the collection of all equivalence classes of weakly continuous
unitary representations of $G$ into $B(\mathcal{H}_{\pi })$ for some Hilbert
space $\mathcal{H}_{\pi }.$

Every $\pi \in \Sigma _{G}$ lifts to a $*$-representation of $L^{1}(G)$ via the formula
	\[\langle \pi (f) \xi ,\eta \rangle =\int_{G} f(x) \langle \pi (x) \xi ,\eta  \rangle dx\]
for any $\xi,\eta\in \Hi_\pi$. We can define a norm on $L^{1}(G)$ by
	\[\| f\| _{*}=\sup \{\| \pi (f)\| : \pi \in \Sigma_G \}.\]
The completion of $L^{1}(G)$ with respect to this norm is the full group $C^*$-algebra, which we denote by $C^{*}(G)$.

A function of the form
	\[u(x)=\la \pi (x) \xi ,\eta \ra (\xi ,\eta \in \mathcal{H})\]
is called a coefficient functions of $\pi.$ Let
	\[ B(G)=\{u(x)= \langle\pi (x)\xi ,\eta \rangle : \pi \in \Sigma_{G},\xi ,\eta \in \mathcal{H}_{\pi }\}. \]
Then $B(G)$ is a commutative algebra of continuous functions on $G$ with respect to pointwise operations.
It is well known that $B(G)$ can be identified with the dual of the group $C^{*}$-algebra $C^{*}(G)$.
In case $G$ is abelian, $C^*(G)$ is the image of $C_0(\widehat {G})$ under the generalized Fourier transform. 
With respect to the dual norm, $B(G)$ is a commutative Banach algebra called the Fourier-Stieltjes algebra of $G.$

The left regular representation $\lambda $ acts on $L^{2}(G)$ as follows:
	\[ (\lambda (y)(f))(x):= f(y^{-1}x) \]
for each $x,y\in G,$ $f\in L^{2}(G).$ We denote by $VN(G)$ the closure of $\text{span}\{\lambda (x) : x\in G\}$ in the weak operator topology of $B(L^{2}(G)).$
$VN(G)$ is a von Neumann algebra called the group von Neumann algebra of $G$. Its predual is $A(G),$ the algebra of continuous functions
that are coefficient functions of the left regular representation $\lambda$ of $G$.
$A(G)$, the Fourier algebra of $G,$ is a closed ideal in $B(G)$
and the norm induced on $A(G)$ as the predual of $VN(G)$ agrees with the norm it inherits from $B(G)$ (\cite{Eymard}).
For any $f\in L^1(G)$ we sometimes denote $\lambda(f)$ by $L_f$ to emphasize that it is the left convolution with respect to $f$ on $L^2(G)$.

A continuous function $f:G\rightarrow \mathbb{C}$ is called a cb-multiplier of $A(G)$
when the multiplication operator with respect to $f$ on $A(G)$ is a completely bounded map.
The collection of all cb-multipliers of $A(G)$ will be denoted by $M_{cb}A(G)$, and $A(G) \subseteq M_{cb}A(G)$.
The closure of $A(G)$ in $M_{cb}A(G)$ will be denoted by $A_{cb}(G)$, and it is a completely contractive Banach algebra.

It will be important for our purposes to note that while $L^{1}(G)$ always has a bounded
approximate identity (BAI), Leptin (\cite{Lep}) showed that $A(G)$ has a BAI if and only if $G$ is amenable.
Moreover, $A_{cb}(G)$ has a BAI if and only if $G$ is weakly amenable (\cite{FRS07}).

Let $C^*_r(G)$ denoted the norm closure of $\{ L_f : f\in L^1(G)\}$ in the space
$B(L^{2}(G))$. $C^*_r(G)$ is called the reduced $C^{*}$-algebra of $G$. When $G$ is amenable,
$C^*(G)$ and $C^{*}_r(G)$ agree, but if $G$ is not amenable, then $C^{*}_r(G)$ is a proper quotient of
$C^*(G)$.

We let $C^*_\delta (G)$ denote the $C^*$-algebra in $B(L^{2}(G))$ generated by
$\{\lambda (x) : x\in G\}$. If $G$ is discrete the clearly $C^*_\delta (G)=C^*_r(G)$.
However for non-discrete groups these algebras are different. When $G$ is abelian,
$C^*_\delta (G)$ is the image of $AP(\widehat {G})$, the algebra of almost periodic functions
on $\widehat {G}$.

We will need to know that if $H$ is an open subgroup of $G$, then there is a natural
completely isometric injection $i:C^*(H)\rightarrow C^*(G)$ (see \cite[Section 5]{BKLS}).

\[UCB(\widehat {G})=\overline{\text{span}}\{u\cdot T : u\in A(G), T\in VN(G)\}\subseteq VN(G).\]
Then $UCB(\widehat {G})$ is called the space of \textit{uniformly continuous functionals on}
$A(G)$.

\subsection{Operator spaces} We will now briefly remind the reader about the basic properties of operator spaces.
We refer the reader to \cite{ER00} for further details concerning the notions presented below.

Let $\mathcal{H}$ be a Hilbert space. Then there is a natural identification between the space
$M_n(B(\mathcal{H}))$ of $n\times n$ matrices with entries in $B(\mathcal{H})$ and the space
$B(\mathcal{H}^n)$. This allows us to define a sequence of norms $\{\| \cdot \| _n\}$ on the spaces
$\{M_n(B(\mathcal{H}))\}$. If $V$ is any subspace of $B(\mathcal{H})$, then the spaces
$M_n(V)$ also inherit the above norm. A subspace $V\subseteq B(\mathcal{H})$ together with the family
$\{\| \cdot \| _n\}$ of norms on $\{M_n(V)\}$ is called a \textit{concrete operator space}. This
leads us to the following abstract definition of an operator space:

\begin{defn}
An operator space is a vector space $V$ together with a family $\{\norm{\cdot}_{n}\}$ of Banach space norms on $M_{n}(V)$ such that for each
$A\in M_{n}(V),B\in M_{m}(V)$ and $[a_{ij}],[b_{ij}]\in M_{n}(\mathbb{C})$

\[
\begin{array}{ll}
i) & \norm{ \left[
\begin{array}{ll}
A & 0 \\
0 & B
\end{array}
\right] }_{n+m} = \max \{\norm{A}_n, \norm{B}_m\} \\ & \\
ii) & \norm{ [a_{ij}]A[b_{ij}]}_n \le \norm{[a_{ij}]} \norm{A}_n \norm{[b_{ij}]}
\end{array}
\]

Let $V,W$ be operator space, $\varphi :V\rightarrow W$ be linear. Then
\[
\norm{ \varphi }_{cb} = \sup_n \{\norm{ \varphi _{n}} \}
\]
where $\varphi _{n}:M_{n}(V)\rightarrow M_{n}(W)$ is given by
\[
\varphi _{n}([v_{ij}])=[\varphi (v_{ij})].
\]

We say that $\varphi $ is completely bounded if $\parallel \varphi \parallel_{cb}<\infty ;$
is completely contractive if $\parallel \varphi \parallel_{cb}\leq 1$ and is a complete isometry if each $\varphi _{n}$ is an isometry.

Given two operator spaces $V$ and $W$, we let $CB(V,W)$ denote the space of all completely
bounded maps from $V$ to $W$. Then $CB(V,W)$ becomes a Banach space with respect to the norm
$\| \cdot \| _{cb}$ and is in fact an operator space via the identification $M_n(CB(V,W))\cong
CB(V,M_n(W))$. We simply write $CB(W)$ for $CB(W,W)$
\end{defn}

It is well-known that every Banach space can be given an operator
space structure, though not necessarily in a unique way. It is also
clear that any subspace of an operator space is also an operator space with respect to the inherited norms.
Moreover, for the duals and the preduals of operator spaces, there are canonical operator space structures.
As such the predual of a von Neumann algebra
and the dual of a $C^*$-algebras respectively, the Fourier and Fourier-Stieltjes algebras
inherit natural operator space structures.

Given two Banach spaces $V$ and $W$, there are many ways to define a norm on the
algebraic tensor product $V\otimes W$. Distinguished amongst such norms is the \textit{Banach space
projective tensor product norm} which we denote by $V\otimes ^\gamma W$. A fundamental property of
the projective tensor product is that there is a natural isometry between $(V\otimes ^\gamma W)^*$ and
$B(V,W^*)$. Given two operator spaces $V$ and $W$, there is an operator space analogue of the projective
tensor product norm which we denote by $V\widehat{\otimes}W$. In this case, we have a natural
complete isometry between $(V\widehat{\otimes}W)^*$ and $CB(V,W^*)$.

The operator space analogue of the injective tensor product will be denoted by $V\otimes _{\min}W$.

\begin{defn}

A Banach algebra $A$ that is also an operator space is called a
\textit{completely contractive Banach algebra} if the multiplication map
\[m:A\widehat{\otimes}A \rightarrow A,\; u\otimes v \mapsto uv\]
is completely contractive. In particular, both $B(G)$ and $A(G)$ are
completely contractive Banach algebras (see \cite{Eymard}).

Let $A$ be a completely contractive Banach algebra. An operator space $X$ is called a
\textit{completely bounded left $A$-module}, if $X$ is a left $A$-module and if
\[\pi _{X} :A\widehat{\otimes} X\rightarrow X\]

with
\[\pi _{X}(u\otimes x)=u\cdot x,\]
is completely bounded.

We say that $X$ is \textit{essential} if $A\cdot X$ is dense in $X$.

We will let the collection of all completely bounded left $A$-modules be denoted
by $A$-\textbf{mod}.

If $X,Y\in A$-\textbf{mod}, then we let $_{A}CB(X,Y)$ denote the space of all
completely bounded left $A$-module maps from $X$ to $Y$.

We can define a \textit{completely bounded right $A$-module}
and a \textit{completely bounded $A$-bimodule} analogously. The collection of all such modules
 will be denoted by \textbf{mod}-$A$ and $A$-\textbf{mod}-$\cdot A$ respectively.

\end{defn}

In general, if $A$ is a completely contractive Banach algebra, then its dual
space $A^*$ is a completely bounded left $A$-module via the action
\[(u\cdot T)(v)=T(vu)\]
for every $u,v\in A$ and $T\in A^*$. Moreover, every closed left
$A$-submodule $Y$ of $A^*$ is also a completely bounded left $A$-module.

Finally, we will need the following operator space analogue of Grothendieck's classical approximation property.
The definition of the operator approximation property is due to Effros and Ruan (\cite{ER00}). The definition
we give is not the original statement of the property but rather has been established as an
equivalent formulation in the same reference.

\begin{defn}\label{def-OAP}
We say that an operator space $X$ has the \textit{Operator Approximation Property} (OAP) if
the natural map
\[J:V^*\widehat{\otimes}\, V\rightarrow V^*\otimes _{\min}V\]
is one-to-one.
\end{defn}

\subsection{$L^p$-spaces and operator space structure}\label{subsec-LpOSS}

Whenever we deal with $L^p$-spaces in this paper we will assume that the reader is familiar with standard materials about
complex interpolation of Banach spaces (\cite{BL76}) and operator spaces (\cite{P96}).

For any von Neumann algebra $\M$ we can define $L^p(\M)$ $(1<p<\infty)$ in the sense of Haagerup.
Kosaki, Terp, and Izumi proved that the Bnach space interpolation $[\M, \M_*]_{\frac{1}{p}}$ is isometric to $L^p(\M)$. 
Using Pisier's complex interpolation theory for operator spaces we can endow an operator space structure on $L^p(\M)$ by
	$${\mathcal O}L^p(\M) = [\M, \M^{op}_*]_{\frac{1}{p}}\;\, \text{(operator space sense)},$$
where $E^{op}$ implies the {\it opposite} of an operator space $E$ (\cite[section 2.10]{P03}).
Note that ${\mathcal O}L^2(\M)$ has Pisier's operator Hilbert space (shortly $OH$) structure. 
We are especially interested in $\M=VN(G)$. However, the usual operator space structure on $A(G)$ is obtained by considering it as the predual of $VN(G)$.
Because of this disagreement we consider the following operator space structure on $L^p(VN(G))$.
	$$L^p(VN(G)) = \left\{ \begin{array}{l} \mathcal{O}L^p(VN(G))^{op} \;\, \text{for}\;\, 1\le p \le 2 \\
	\mathcal{O}L^p(VN(G)) \;\, \text{for}\;\, 2\le p \le \infty. \end{array} \right.$$
If we use reiteration theorem, then we get
	$$\mathcal{O}L^p(VN(G)) = [L^2(VN(G))_{oh}, VN(G)^{op}_*]_{\frac{2}{p}-1}$$
for $1\le p\le 2$, where $H_{oh}$ is the operator Hilbert space defined on a Hilbert space $H$.
Since the opposite of $OH$ is still $OH$, we get
	$$L^p(VN(G)) = [L^2(VN(G))_{oh}, A(G)]_{\frac{2}{p}-1}\; \text{for}\; 1\le p\le 2.$$

We are also intersted in $L^p(G)$, the $L^p$-spaces associated to the commutative von Neumann algebra $L^\infty(G)$.
In this case we have $(L^\infty(G))^{op} = L^\infty(G)$, thus taking opposite does not affect to the operator space structure.

\subsection{Operator Segal algebras}

Segal algebras were first defined by H. Reiter for group algebras;
see \cite{RS}, for example. The definition of operator Segal
algebras appeared in \cite{FSW}. However, our abstract definition
deviates from the one given in \cite{FSW} in the sense that we
demand that Segal algebras be essential modules.

Let $A$ be a completely contractive Banach algebra.  An {\it operator
Segal algebra} is a subspace $B$ of $A$ such that

(i) $B$ is dense in $A$,

(ii) $B$ is a left ideal in $A$, and

(iii) \parbox[t]{4.4in}{$B$ admits an operator space structure
$\|\cdot\|_B$ under which it is complete and a completely
contractive $A$-module.}

(iv) $B$ is an essential $A$-module:  $A\cdot B$ is $\|\cdot\|_B$-dense
in $B$.

We will discuss two specific types of operator Segal algebras in the
Fourier algebra $A(G)$.  One is the Lebesgue-Fourier algebra, $S^1A(G) = L^1(G)\cap A(G)$,
whose study was initiated in \cite{GL1}, and which was shown to be
an operator Segal algebra in \cite{FSW}.  The second is
Feichtinger's algebra $S_0(G)$, whose study in the non-commutaive
case was taken up in \cite{Spr2}; this study included an exposition
of the operator space structure.  Though slightly different
terminology was used in that article, it was proved there that
$S_0(G)$ is an operator Segal algebra in $A(G)$, in the sense
defined above. It is important to realize that both $S^1A(G)$ and
$S_0(G)$ are also operator Segal algebras in the group algebra $L^1(G)$
with the convolution product.

On the other hand, $A(G)$ itself can be understood as an operator Segal algebra of $A_{cb}(G)$.

\subsection{General theory about projectivity}

In this subsection we will establish or recall some basic properties of projective modules that we will need in our study.

\begin{defn}
Let $A$ be a completely contractive Banach algebra. Let $X$ and $Y$ be completely bounded
left $A$-modules. A map $T\in CB(X,Y)$ is \textit{admissible} if $ker T$ is completely
complemented in $X$ and the range of $T$ is closed and completely complemented in $Y$.

Let $X\in A$-{\bf mod}. Then $X$ is said to be \textit{operator projective in}
$A$-{\bf mod} if whenever $E,F\in A$-{\bf mod} and
$T\in {_{A}CB(E,F)}$ is admissible and surjective,
then for each $S\in {_{A}CB(X,F)}$, there exists $R\in {_{A}CB(X,E)}$ such that
$T\circ R=S$.

\begin{equation*}
\begin{tabular}{lll}
&  & $X$ \\
& $\overset{R}{\swarrow }$ & $\downarrow S$ \\
$E$ & $\underset{T}{\longrightarrow }$ & $F$%
\end{tabular}
\end{equation*}

Equivalently, $X$ is operator projective in $A$-{\bf mod} if $Z$ is any
submodule of $F$, then every $T\in {_{A}CB(X,F/Z)}$ lifts to
a map in $T\in {_{A}CB(X,F)}$

We can also define operator projectivity in ${\bf mod}$-$A$ and in $A$-${\bf mod}$-$A$
in a similar manner. In this case, we say that $A$ is \textit{operator biprojective} if $A$ is
operator projective in $A$-${\bf mod}$-$A$.
\end{defn}

Let $A$ be a completely contractive Banach algebra. Then we denote the unitization $A \oplus \mathbb{C}$ of $A$ by $A_+$.
Let $X$ be a completely bounded left $A$-module with the multiplication map
	$$\pi_X : A \prt\, X\rightarrow X.$$
Then $\pi_X$ can be extended to $\pi_{X,+} : A_+ \prt\, X \rightarrow X$ in a canonical way.
Using $\pi_{X,+}$ we have a useful characterization of operator projectivity by P. Woods (\cite[Corollary 3.20, Proposition 3.24]{W02}).

	\begin{prop}\label{prop-woods}
	Let $A$ be a completely contractive Banach algebra and $X$ be a completely bounded left $A$-module.
	Then, $X$ is operator projective if and only if we have a completely bounded $A$-module map
		$$\rho : X \rightarrow A_+ \prt\, X,$$
	which is a right inverse of the extended multiplication map $\pi_{X,+}$.
	
	When $X$ is essential, $X$ is operator projective if and only if we have a completely bounded $A$-module map
		$$\rho : X \rightarrow A \prt\, X,$$
	which is a right inverse of the multiplication map $\pi_X$.
	\end{prop}

The following theorem, which is proved in \cite[Theorem 3.2]{FLS}, gives a sufficient condition for
a module to be operator projective.

	\begin{thm}\label{thm-essential}
	Let $A$ be an operator biprojective completely contractive Banach algebra with a BAI,
	and let $X$ be an essential completely bounded left $A$-module.
	Then $X$ is operator projective in $A$-{\bf mod}.
	\end{thm}

It is important to note that the assumption that the completely contractive
Banach algebra $A$ has a BAI is crucial in the previous two statements.
As we have previously observed, while $L^{1}(G)$ always has a bounded approximate identity,
neither of $S^1A(G)$ nor $S_0(G)$ has a BAI.
Since $A_{cb}(G)$ has a BAI if and only if $G$ is weakly amenable,
we will see that most of the positive results we obtain with respect to identifying
operator projective $A_{cb}(G)$-modules will require the assumption of weak amenability.

We finish this section with the following proposition which may be viewed as an operator space analogue of \cite[Corollary 4.5]{Hel}.
Its proof is presented in \cite[Proposition 3.5]{FLS}.

	\begin{prop}\label{prop-proj}
	Let $A$ be a completely contractive Banach algebra, and let $X$ be a completely bounded left $A$-module.
	Suppose that $X$ is operator projective in $A$-{\bf mod} and $X$ or $A$ have OAP.
	Then for any non-zero element $x\in X$, there is a map $T \in {}_A CB(X, A_+)$ such that $T(x) \neq 0$.
	\end{prop}

\subsection{$A(G)$-module structures on various spaces}
In this subsection we recall $A(G)$-module structures on various spaces (see \cite{FLS} for the details),
which will be extended to the case of other algebras later in this paper. Note that all the module operations below are completely contractive.

First of all, $A(G)$ is a $A(G)$-module by the algebra operation
	$$\pi_1 : A(G)\prt A(G)\rightarrow A(G),\; f \otimes g \mapsto fg,$$
and $A(G)^{**}$ have the natural $A(G)$-module structures as the bidual.
$VN(G)$ has the dual module structure given by
	$$\pi_\infty : A(G)\prt VN(G)\rightarrow VN(G),\; f \otimes L_g \mapsto L_{fg}.$$
$C^*_r(G)$ and $UCB(\widehat{G})$ are understood as submodules of $VN(G)$. Similarly, $C^*_\delta(G)$ is a submodule of $VN(G)$
	$$\pi_\delta : A(G)\prt C^*_\delta(G) \rightarrow C^*_\delta(G),\; f \otimes \lambda(x) \mapsto f(x)\lambda(x).$$
Let $i : L^1(G) \hookrightarrow C^*(G)$ is the canonical embedding induced from the universal representation of $G$.
Then $\{i(f) : f\in L^1(G)\}$ is dense in $C^*(G)$.
Now we can consider a natural $A(G)$-module structure on $C^*(G)$ by
	$$\pi_{C^*(G)} : A(G) \prt C^*(G) \rightarrow C^*(G),\; f \otimes i(g) \mapsto i(f g).$$
$UCB(\widehat{G})^*$ has the dual module structure.

The $A(G)$-module structure on $L^p(G)$ $(1\le p \le \infty)$ is given by
	\begin{equation}\label{module-Lp}
	\pi_{L^p(G)} : A(G) \prt L^p(G) \rightarrow L^p(G),\; f\otimes g \mapsto f g.
	\end{equation}
It would be beneficial for us later to explain how we get $\pi_{L^p(G)}$.
Indeed, since $L^\infty(G)$ is a completely contractive Banach algebra under the pointwise multiplication, we have a complete contraction
	$$\Pi_\infty : L^\infty(G) \prt L^\infty(G) \rightarrow L^\infty(G),\; f\otimes g \mapsto f g.$$
On the other hand $L^1(G)$ is the predual of $L^\infty(G)$,
so that it is a $L^\infty(G)$-module with the following completely contractive map
	$$\Pi_1 : L^\infty(G) \prt L^1(G) \rightarrow L^1(G),\; f\otimes g \mapsto f g.$$
Then, we get a complete contraction
	\begin{equation}\label{com-Lp}
	\Pi_p : L^\infty(G) \prt L^p(G) \rightarrow L^p(G),\; f\otimes g \mapsto f g
	\end{equation}
by the bilinear complex interpolation (\cite[section 2.7]{P03}).
Since the formal identity $j : A(G) \hookrightarrow L^\infty(G)$ is a complete contraction we get $\pi_{L^p(G)} = (j\otimes id_{L^p(G)})\circ \Pi_p$.

Finally, we consider the $A(G)$-module structure on $L^p(VN(G))$ $(1< p <\infty)$.
Let $\varphi$ be the Plancherel weight on $VN(G)$ (\cite[III.3.3]{Bla} or \cite[VII.3]{Ta03}) given by
	$$\varphi(L_f) = f(e)\;\, \text{for any}\;\, f\in C_c(G).$$
Let $n_\varphi = \{x\in \M : \varphi(x^*x) <\infty\}$ and $\Lambda$ be the canonical embedding of $n_\varphi$ into $L^2(VN(G))$,
where $L^2(VN(G))$ is the Hilbert space obtained by the completion of $n_\varphi$ with the inner product
	$\left\langle x, y\right\rangle_\varphi = \varphi(y^* x)$ for any $x, y \in n_\varphi$.
We will use complex interpolation again, but unfortunately the maps $\pi_1$ and $\pi_\infty$
are not compatible in the sense of interpolation theory (see \cite[section 6]{FLS} for the details).
Instead, we consider the following two intermidiate maps
	$$\pi_2 : A(G) \prt L^2(VN(G)) \rightarrow L^2(VN(G)),\; f \otimes \Lambda(L_{\check{g}}) \mapsto \Lambda(L_{\check{f} \check{g}})$$
and
	$$\pi'_2 : A(G) \prt L^2(VN(G)) \rightarrow L^2(VN(G)),\; f \otimes \Lambda(L_g) \mapsto \Lambda(L_{f g})$$
where $\check{f}$ is given by $\check{f}(x) = f(x^{-1})$ for any $x\in G$.
If we recall the complete isometry
	$$L^2(G) \rightarrow L^2(VN(G)),\; f\mapsto L_f$$
for $f\in C_c(G)$, it is easy to see that $\pi'_2$ is a complete contraction by \eqref{module-Lp}.
The case of $\pi_2$ can be shown similarly if we repeat the above argument for $\Pi_p$
with the $L^\infty(G)$-module structure on $L^\infty(G)$ given by the complete contraction
	$$L^\infty(G) \prt L^\infty(G) \rightarrow L^\infty(G),\; f\otimes \check{g} \mapsto \check{f}\check{g}$$
instead of $\Pi_\infty$. Since $\pi_2$ and $\pi'_2$ are compatible with $\pi_1$ and $\pi_\infty$, respectively,
we get the $A(G)$-module actions on $L^p(VN(G))$ $(1<p<2)$ by the interpolation of $\pi_1$ and $\pi_2$
and on $L^p(VN(G))$ $(2\le p<\infty)$ by the interpolation of $\pi_\infty$ and $\pi'_2$, respectively.

Note that the above explanation for the complete contractivities of $\pi_2$ and $\pi'_2$ is much simpler than that of \cite[section 6]{FLS},
but this approach can not be extended to the case of Kac algebras and does not cover the $A(G)$-module structure on $L^2(G)_c$ (resp. $L^2(G)_r$),
where $H_c$ (resp. $H_r$) implies the column (resp. row) Hilbert space on $H$ for any Hilbert space $H$.

\section{Operator approximation property of Segal algebras}

As it is pointed out in Proposition \ref{prop-proj}, it is beneficial if we can determine
whether a given algebra or a module has operator approximation property
(OAP). In \cite{Z1} and \cite{Z2}, Y. Zhang introduced the concept of being approximately
complemented for a subspace $E$ of a normed space $X$ and showed that it is closely
related to $E$ having the approximation property. In this section, we modify this concept,
so that it can be viewed in the category of operator spaces and we apply it to deduce OAP for certain Segal algebras.

Let $\K$ denote the algebra of compact operators on $\ell^2(\mathbb{N})$. We recall from \cite{ER00} that,
for an operator space $V$, $\K(V):=\K \otimes_{\min} V$ denote the operator injective tensor
product of $\K$ and $V$.

	\begin{defn}
	Let $(V, \norm{\cdot}_{M_n(V)})$ be an operator space,
	and let $W$ be a linear subspace of $V$ with a (possibly different) operator space structure $\{\norm{\cdot}_{M_n(W)}\}$
	for which the inclusion $\iota : W \hookrightarrow V$ is completely bounded.
	We say that $W$ is {\it completely approximately complemented} in $V$ if there is a net $\{\varphi_i\} \in CB(V,W)$ such that
	for every $w\in \K(W)$,
		$$\norm{\varphi_{i,\infty} (w)-w} \stackrel{i}{\longrightarrow} 0,$$
	where $\varphi_{i,\infty} = id_\K \otimes \varphi_i$.
	\end{defn}
Note that $\K(W)$ can be viewed as a subspace of $\K(V)$.
It is easy to observe that $W$ is completely approximately complemented in $V$ if and only if
for every $\epsilon >0$ and $w_1,\ldots, w_n \in \K(W)$, there is $\varphi\in CB(V,W)$ such that
	\begin{equation}\label{cb-app-complemented}
	\norm{\varphi_\infty(w_k)-w_k} <\epsilon,\; 1\le k\le n,
	\end{equation}
where $\varphi_\infty = id_\K \otimes \varphi$.
Since $\K\oplus \cdots \oplus\K$ can be regarded as a subspace of $\K$,
it suffices to assume that \eqref{cb-app-complemented} holds for a single $w\in \K(W)$.

	\begin{thm}\label{thm-CAC-OAP}
	Let $V$ be an operator space with OAP. If $W$ is completely approximately complemented in $V$, then $W$ has OAP.
	\end{thm}
\begin{proof}
Suppose that $W$ is completely approximately complemented in $V$. Let $w\in \K(W)$ and $\epsilon > 0$.
Then there is $\varphi \in CB(V,W)$ such that
	$$\norm{\varphi_\infty(w) -w}_{\K(W)}<\epsilon.$$
Since $V$ has OAP, there is $\psi\in CB(V,V)$ with finite rank such that
	$$\norm{\psi_\infty(w)-w}_{\K(V)} < \frac{\epsilon}{\norm{\varphi}_{cb}}.$$
Put $\phi = \varphi\circ \psi \circ \iota$, where $\iota:W\hookrightarrow V$ is the inclusion.
Then $\phi \in CB(W,W)$ with finite rank and
	\begin{align*}
	\norm{\phi_\infty(w)-w}_{\K(W)} & \le \norm{\varphi_\infty(\psi_\infty(w)-w)}_{\K(W)}
	+ \norm{\varphi_\infty(w)-w}_{\K(W)}\\
	& \le \norm{\varphi_\infty}\norm{\psi_\infty(w)-w}_{\K(V)} + \epsilon \le 2\epsilon.
	\end{align*}
Thus $W$ has OAP.
\end{proof}

	\begin{thm}\label{thm-OAP-Segal}
	Let $A$ be a completely contractive Banach algebra with a BAI, and let $B$ be an operator Segal algebra of $A$. Then
		\begin{enumerate}
			\item $B$ is completely approximately complemented in $A$.
			\item $A$ has OAP if and only if $B$ has OAP.
		\end{enumerate}
	\end{thm}
\begin{proof}
(1) Let $\{e_\alpha\}_{\alpha \in I}$ be a BAI in $A$. Since $B$ is essential, $\{e_\alpha\}_{\alpha \in I}$ is also a BAI for $B$ in $A$.
That is
	$$\norm{e_\alpha b - b}_B \rightarrow 0\;\, \text{as}\;\, \alpha \rightarrow \infty$$
for any $b\in B$.
Now for each $\alpha \in I$, we define a completely bounded map
	\begin{equation}\label{BAI-maps}
	\varphi_\alpha : A \rightarrow B,\; a \mapsto e_\alpha a.
	\end{equation}
Let $M= \sup\{\norm{e_\alpha}_A : \alpha\in I\}$ and $\epsilon > 0$ be arbitrary.
Then, for any $\widetilde{b} = (b_{ij}) \in \K(B)$, there is $n\in \n$ such that
	$$\norm{\widetilde{b}-\widetilde{b}_n}_{\K(B)} < \min\{\epsilon, \frac{\epsilon}{M}\},$$
where $\widetilde{b}_n = (b_{ij})^n_{i,j=1}$.

Hence, for every $\alpha\in I$,
	\begin{align*}
	\norm{\varphi_{\alpha,\infty}(\widetilde{b}) - \widetilde{b}}_{\K(B)}
	& \le \norm{\varphi_{\alpha,\infty}(\widetilde{b}-\widetilde{b}_n)}
	+ \norm{\varphi_{\alpha,\infty}(\widetilde{b}_n)-\widetilde{b}_n} + \norm{\widetilde{b}_n-\widetilde{b}}\\
	& \le \norm{\varphi_{\alpha,\infty}}_{cb}\norm{\widetilde{b}-\widetilde{b}_n} + \epsilon
	+ \norm{(e_\alpha b_{ij}-b_{ij})^n_{i,j=1}}_{\K(B)}\\
	& \le 2\epsilon + \sum^n_{i,j=1}\norm{e_\alpha b_{ij}-b_{ij}}_B.
	\end{align*}
Thus there is $\alpha_0 \in I$ such that $\alpha \ge \alpha_0$ implies
	$$\norm{\varphi_{\alpha,\infty}(\widetilde{b}) - \widetilde{b}}_{\K(B)} < 3\epsilon.$$

\vspace{0.3cm}
(2) ($\Rightarrow$) It is clear by Theorem \ref{thm-CAC-OAP} and part (1).

($\Leftarrow$) Let $\epsilon >0$ and $\widetilde{a} = (a_{ij})\in \K(A)$.
Let $\{e_\alpha\}_{\alpha\in I}$ be a BAI in $A$ and recall the maps $\varphi_\alpha$'s in \eqref{BAI-maps}.
Then we can show that there is $\beta \in I$ such that
	$$\norm{\varphi_{\beta, \infty}(\widetilde{a})-\widetilde{a}}_{\K(A)} < \epsilon$$
as in the proof of part (1). Since $B$ has OAP, there is $\psi\in CB(B,B)$ with finite rank such that
	$$\norm{\psi_\infty(\varphi_{\beta, \infty}( \widetilde{a})) - \varphi_{\beta, \infty}( \widetilde{a})}_{\K(B)}<\epsilon.$$
Now we set $\phi = \psi \circ \varphi_{\beta} \in CB(A,B) \subseteq CB(A,A)$.
Clearly $\phi$ has finite rank. Moreover,
	\begin{align*}
	\norm{\phi_\infty(\widetilde{a})-\widetilde{a}}_{\K(A)}
	& = \norm{\psi_\infty(\varphi_{\beta, \infty}(\widetilde{a}))-\widetilde{a}}_{\K(A)}\\
	& \le \norm{\psi_\infty(\varphi_{\beta, \infty}(\widetilde{a}))-\varphi_{\beta, \infty}(\widetilde{a})}_{\K(A)}
	+ \norm{\varphi_{\beta, \infty}(\widetilde{a}) - \widetilde{a}}_{\K(A)}\\
	& \le \norm{\psi_\infty(\varphi_{\beta, \infty}(\widetilde{a}))-\varphi_{\beta, \infty}(\widetilde{a})}_{\K(B)} + \epsilon < 2\epsilon.
	\end{align*}
Thus $A$ has OAP.
\end{proof}

Since $L^1(G)$ always has OAP and a BAI, we get the following corollary.

	\begin{cor}
	Let $G$ be a locally compact group. Then every operator Segal algebra of $L^1(G)$ has OAP.
	In particular, $S^1A(G)$ and $S_0(G)$ have OAP.
	\end{cor}

	\begin{cor}\label{cor-OAP-Acb(G)}
	Let $G$ be a weakly amenable locally compact group. Then $A_{cb}(G)$ has OAP.
	\end{cor}
\begin{proof}
Since $G$ is weakly amenable, $A_{cb}(G)$ has a BAI (\cite{FRS07}), and $A(G)$ has OAP (\cite{KR}).
Also, $A(G)$ is an operator Segal algebra of $A_{cb}(G)$. Thus the result follows from Theorem \ref{thm-OAP-Segal}.
\end{proof}

\section{operator Segal algebras of the Fourier algebra}\label{S:operator-SA(G)-S_0(G)}

In this section, we study operator left projectivity of various modules over
operator Segal algebras of the Fourier algebra. We will specially focus
on the following two algebras: $S^1A(G) = A(G) \cap L^1(G)$, the Lebesgue-Fourier algebra,
and $S_0(G)$, the Feichtinger's Segal algebra.

Let $S(G)$ be an operator Segal algebra of $A(G)$ with the inclusion $j : S(G) \rightarrow A(G)$.
If $X$ is a completely bounded left $A(G)$-module with the multiplication map
	$$\pi : A(G) \prt X \rightarrow X,$$
then we can impose a natural $S(G)$-module structure on $X$ by the multiplication map
	$$\widetilde{\pi} : S(G) \prt X \stackrel{j \otimes id_X}{\longrightarrow} A(G) \prt X \stackrel{\pi}{\longrightarrow} X.$$

	\begin{prop}\label{P:op. S(G)-A(G)}
	Let $S(G)$ be an operator Segal algebra of $A(G)$.
	If $X \in$ $A(G)$-{\bf mod} and operator projective in $S(G)$-{\bf mod}, then $X$ is operator projective in $A(G)$-{\bf mod}.
	\end{prop}

\begin{proof}
Suppose that $X$ is operator projective in $S(G)$-\m. Hence there is a completely
bounded left $S(G)$-module morphism $\rho : X \to S(G)_+ \prt X$ such that $\widetilde{\pi}_+ \circ
\rho=id_X$. Then $\pi_+ \circ (j_+ \otimes id_X) \circ \rho=id_X$. Moreover $(j_+ \otimes id_X) \circ \rho$
is an $A(G)$-module map since $S(G)$ is dense in $A(G)$. Thus $X$ is operator projective in $A(G)$-\m.
\end{proof}

As it is shown in \cite[Table 1]{FLS}, most of the modules that we are interested
to study fail to be operator projective in $A(G)$-\m\ for a non-discrete $G$. By comparing this fact
with the preceding proposition, it follows that the same will hold for operator projectivity in
$S(G)$-\m. Therefore, in order to obtain the complete characterization, we mainly have to consider the case when $G$ is discrete.
On the other hand, for a discrete group $G$, it is routine to verify that $S(G)$
will have the algebra $\ell^1(G)$, with the pointwise product, as an operator Segal algebra. Moreover, in most
of the desirable cases such as $S^1A(G)$ and $S_0(G)$, $S(G)$ is actually $\ell^1(G)$ for discrete groups (\cite[Corollary 2.5]{Spr2}).
Thus we are reduced to the study of operator projectivity in $\ell^1(G)$-\m.
We will show that operator projectivity in $\ell^1(G)$-\m\ rarely happens unless $G$ is finite.

The following map will be used frequently in this section.
Note that $\ell^1(G)\prt X$ can be identified with the vector valued $\ell^1$-space $\ell^1(G;X)$ for any Banach space $X$.

	\begin{prop} Let
		\begin{equation}\label{projection}
		\Theta : \ell^1(G;X) \rightarrow \ell^1(G;X),\; u\mapsto \Theta(u)
		\end{equation}
	with $\Theta(u)(s) = \delta_s \cdot u(s)$ for $s\in G$. Then $\Theta$ is a linear (complete)
    contraction.
	\end{prop}

\begin{proof}
For any $u \in \ell^1(G;X)$ we have
	$$\norm{\Theta(u)} = \sum_{s\in G}\norm{\delta_s \cdot u(s)}_X \le \sum_{s\in G}\big[\norm{\delta_s}_{\ell^1(G)}\norm{u(s)}_X\big] = \norm{u}_{\ell^1(G;X)}.$$
\end{proof}

\subsection{The modules $C^*_r(G)$, $UCB(\widehat{G})$, $C^*_\delta(G)$, $VN(G)$, and $C^*(G)$}

In this section, we show that the operator projectivity of $C^*_r(G)$, $UCB(\widehat{G})$, $VN(G)$,
$C^*(G)$, or $C^*_\delta(G)$ in $S^1A(G)$-\m\ or $S_0(G)$-\m\ implies that $G$ is finite. We start with the discrete case.

	\begin{thm}\label{thm-discrete-Segal-VN(G)}
	Let $G$ be a discrete group.
		\begin{enumerate}
			\item Let $X$ be a completely bounded $\ell^1(G)$-submodule of $VN(G)$	containing $\{\lambda(s) : s\in G\}$.
			Then $X$ is operator projective in $\ell^1(G)$-{\bf mod} if and only if $G$ is finite.
			
			\item $C^*(G)$ is operator projective in $\ell^1(G)$-{\bf mod} if and only if $G$ is finite.
		\end{enumerate}
	\end{thm}
\begin{proof}
(1) When $G$ is finite, the operator projectivity of $X$ is trivial.
Now suppose that $X$ is operator projective in $\ell^1(G)$-{\bf mod}. There is a completely bounded left $\ell^1(G)$-module map
	$$\rho: X \rightarrow \ell^1(G)_+ \prt X$$
such that $\pi_+ \circ \rho = id_X$.
By a standard argument (see, for example, the proof of \cite[Theorem 5.9]{FLS}) we can show that,
for every $s\in G$,
	$$\rho(\lambda(s)) = \delta_s\otimes (\lambda(s) + x_s)$$
for some $x_s \in Y$, where $X = Y\oplus \mathbb{C}\lambda(s)$.
Now we set
	$$\Psi = \Theta \circ \rho : X_1 \rightarrow \ell^1(G;X),$$
where $X_1$ is the closure of $\{\lambda(s) : s\in G\}$ in $X\subseteq VN(G)$. Then we have
	$$\Psi(\lambda(s)) = \delta_s \otimes \lambda(s)$$
for any $s\in G$ since $\delta_s \cdot x_s = 0.$
Thus, for any $f\in \ell^1(G)$
	$$\Psi(\lambda(f)) = \sum_{s\in G}f(s)\delta_s \otimes \lambda(s),$$
and
	\begin{align*}
	\norm{\lambda(f)}_X & \le \norm{f}_1 = \sum_{s\in G}\abs{f(s)} = \sum_{s\in G} \norm{f(s)\lambda(s)}_X = \norm{\Psi(\lambda(f))}_{\ell^1(G;X)}\\
	& \le \norm{\Psi}\norm{\lambda(f)}_X.
	\end{align*}
Consequently, $\norm{\cdot}_1$ and $\norm{\cdot}'$ are equivalent on $\ell^1(G)$, where $\norm{f}' = \norm{\lambda(f)}_{VN(G)}$ for $f\in \ell^1(G)$.
In particular, $(\ell^1(G), \norm{\cdot}')$ is complete, and so it is a $C^*$-algebra with respect to the convolution product.
Thus $\ell^1(G)$ with convolution is Arens regular, which implies that $G$ is finite.

\vspace{0.3cm}
(2) Let $i : \ell^1(G) \rightarrow C^*(G)$ be the canonical embedding.
Using a similar argument as above we can show that $\norm{\cdot}_1$ and $\norm{\cdot}''$ are equivalent on $\ell^1(G)$,
where $\norm{f}'' = \norm{i(f)}_{C^*(G)}$ for $f\in \ell^1(G)$.
In particular, $(\ell^1(G), \norm{\cdot}'')$ is complete, and so it is a $C^*$-algebra with respect to the convolution product.
Thus $\ell^1(G)$ with convolution is Arens regular again, which implies that $G$ is finite.

\end{proof}

	\begin{thm}
	Let $A = S^1A(G)$ or $S_0(G)$.
	Then $C^*_r(G)$ (resp. $UCB(\widehat{G})$, $VN(G)$, $C^*(G)$, or $C^*_\delta(G)$) is operator projective in $A$-{\bf mod} if and only if $G$ is finite.
	\end{thm}

\begin{proof}
$``\Leftarrow"$ Clear.\\
$``\Rightarrow"$ It follows from Proposition \ref{P:op. S(G)-A(G)} and \cite[Theorems 5.4, 7.7, and 7.10]{FLS}
that $G$ must be discrete. Therefore $G$ is finite by Theorem \ref{thm-discrete-Segal-VN(G)}.
\end{proof}

\subsection{The modules $A(G)$, $A(G)^{**}$, $UCB(\widehat{G})^*$, and $L^p(VN(G))$ ($1<p<\infty$)}

Similar to the preceding section, we show that operator projectivity of $A(G)^{**}$, $UCB(\widehat{G})^*$,
and $L^p(VN(G))$ ($2\leq p<\infty$) in $S^1A(G)$-\m\ or $S_0(G)$-\m\ implies that $G$ is finite.

\begin{thm}\label{T:disceret-L(VN(G))}
	Let $G$ be a discrete group.
	Then $L^p(VN(G))$ $(1\le p<\infty)$ (resp. $A(G)^{**}$ and $UCB(\widehat{G})^*$) is operator projective in $\ell^1(G)$-{\bf mod} if and only if $G$ is finite.
	\end{thm}
\begin{proof}
We basically follow the proof of Theorem \ref{thm-discrete-Segal-VN(G)}.
We consider the case $2\le p <\infty$ (the case $1\le p <2$ is similar).
Since $G$ is discrete, $VN(G)$ is a finite von Neumann algebra equipped with the canonical trace $\tau$
satisfying $\tau(1_{VN(G)}) = 1$. Moreover,
	$$VN(G) \subseteq L^p(VN(G))$$
contractively for $1\le p <\infty$.

We recall the contraction $\Theta$ in \eqref{projection} and suppose that $L^p(VN(G))$ is operator projective in $\ell^1(G)$-{\bf mod}.
Then there is a completely bounded left $\ell^1(G)$-module map
	$$\rho: L^p(VN(G)) \rightarrow \ell^1(G) \prt L^p(VN(G))$$
such that $\pi \circ \rho = id_{L^p(VN(G))}$.

By a standard argument we can show that, for any $s\in G$,
	$$\rho(\lambda(s)) = \delta_s\otimes (\lambda(s) + x_s)$$
for some $x_s \in Y$, where $L^p(VN(G)) = Y\oplus \mathbb{C}\lambda(s)$.
Now we set
	$$\Psi = \Theta \circ \rho : L^p(VN(G)) \rightarrow \ell^1(G;L^p(VN(G))).$$
Then we have for any $f\in \ell^1(G)$,
	$$\Psi(\lambda(f)) = \sum_{s\in G}f(s)\delta_s \otimes \lambda(s),$$
and since $\norm{\lambda(s)}_{L^p(VN(G))} = 1$ for any $s\in G$, we have
	\begin{align*}
	\norm{\lambda(f)}_{VN(G)} & \le \norm{f}_1 = \sum_{s\in G}\abs{f(s)}\\
	& = \sum_{s\in G} \norm{f(s)\lambda(s)}_{L^p(VN(G))} = \norm{\Psi(\lambda(f))}_{\ell^1(G;L^p(VN(G)))}\\
	& \le \norm{\Psi}\norm{\lambda(f)}_{L^p(VN(G))} \le \norm{\Psi}\norm{\lambda(f)}_{VN(G)}.
	\end{align*}
Consequently, as in the proof of Theorem \ref{thm-discrete-Segal-VN(G)}, the norms $\norm{\cdot}_1$ and $\norm{\cdot}'$ are equivalent on $\ell^1(G)$.
Thus for the same reason $G$ is finite.
A similar argument shows that if either of $A(G)^{**}$ or $UCB(\widehat{G})^*$ is operator projective in $\ell^1(G)$-\m,
the $\|\cdot\|_1$ and $\|\cdot\|_{A(G)}$ are equivalent on $\ell^1(G)$. Thus $G$ must be finite.
\end{proof}

\begin{thm}
	Let $A = S^1A(G)$ or $S_0(G)$.
	Then $L^p(VN(G))$ $(2\leq p<\infty)$ $($resp. $A(G)^{**}$ and $UCB(\widehat{G})^*$) is operator projective in $A$-{\bf mod} if and only if $G$ is finite.
	\end{thm}
\begin{proof}
$``\Leftarrow"$ Clear.\\
$``\Rightarrow"$ It follows from Proposition \ref{P:op. S(G)-A(G)} and \cite[Theorems 4.9, 4.10, and 6.9]{FLS} that $G$ must be discrete.
Therefore $G$ is finite by Theorem \ref{T:disceret-L(VN(G))}.
\end{proof}

\subsection{The modules $L^p(G)$ $(1\le p \le \infty)$}

We first note that when $G$ is discrete, $\ell^1(G)$ with pointwise multiplication is operator biprojective. Thus $\ell^1(G)$ is operator projective in $\ell^1(G)$-{\bf mod}.

	\begin{thm}\label{T:l^p-discrete}
	Let $G$ be a discrete group.
	Then $\ell^p(G)$ ($1< p \le \infty$) is operator projective in $\ell^1(G)$-{\bf mod} if and only if $G$ is finite.
	\end{thm}
\begin{proof}
We again follow the proof of Theorem \ref{thm-discrete-Segal-VN(G)}.
We recall the contraction $\Theta$ in \eqref{projection} and suppose that $\ell^p(G)$ is operator projective in $\ell^1(G)$-{\bf mod}.
Then there is a completely bounded left $\ell^1(G)$-module map
	$$\rho: \ell^p(G) \rightarrow \ell^1(G)_+ \prt \ell^p(G)$$
such that $\pi_+ \circ \rho = id_{\ell^p(G)}$.

By a standard argument we can show that, for any $s\in G$,
	$$\rho(\delta_s) = \delta_s\otimes (\delta_s + x_s)$$
for some $x_s \in Y$, where $\ell^p(G) = Y\oplus \mathbb{C}\delta_s$.
Now we set
	$$\Psi = \Theta \circ \rho|_{\ell^1(G)} : \ell^1(G) \rightarrow \ell^1(G;\ell^p(G)).$$
Then we have for any $f\in \ell^1(G)$
	$$\Psi(f) = \sum_{s\in G}f(s) \delta_s \otimes \delta_s,$$
and since $\norm{\delta_s}_{\ell^p(G)} = 1$ for any $s\in G$, we have
	\begin{align*}
	\norm{f}_p & \le \norm{f}_1 = \sum_{s\in G}\abs{f(s)}\\
	& = \sum_{s\in G} \norm{f(s)\delta_s}_{\ell^p(G)} = \norm{\Psi(f)}_{\ell^1(G;\ell^p(G))}\\
	& \le \norm{\Psi}\norm{f}_p.
	\end{align*}
Consequently, $\norm{\cdot}_1$ and $\norm{\cdot}_p$ are equivalent on $\ell^1(G)$, which implies that $G$ is finite.
\end{proof}

\begin{thm}
	Let $A = S^1A(G)$ or $S_0(G)$. Then:\\
 $($i$)$ $L^1(G)$ is operator projective in $A$-{\bf mod} if and only if $G$ is discrete.\\
 $($ii$)$ $L^p(G)$ $(1 < p <\infty)$ is operator projective in $A$-{\bf mod} if and only if $G$ is finite.
	\end{thm}

\begin{proof}
(i) $``\Leftarrow"$ It follows from the fact that $\ell^1(G)$ is operator biprojective.\\
$``\Rightarrow"$ It follows from Proposition \ref{P:op. S(G)-A(G)} and \cite[Theorems 7.1]{FLS}.\\
(ii) $``\Leftarrow"$ Clear.\\
     $``\Rightarrow"$ It follows from Proposition \ref{P:op. S(G)-A(G)} and \cite[Theorems 7.1]{FLS} that $G$ must be discrete. Therefore $G$ is finite by Theorem \ref{T:l^p-discrete}.
\end{proof}

\section{The algebra $A_{cb}(G)$}

In this section, we study the operator projectivity of the modules we considered
in Section \ref{S:operator-SA(G)-S_0(G)} but this time as left modules of $A_{cb}(G)$,
the closure of $A(G)$ in $M_{cb}A(G)$.

\subsection{The modules $C^*_r(G)$, $UCB(\widehat{G})$, $C^*_\delta(G)$, $VN(G)$, and $C^*(G)$}

We first need to explain the completely bounded $A_{cb}(G)$-module structures on each
one of these spaces.
Since $A_{cb}(G)$ is the closure of $A(G)$ in $CB(A(G))$ and we have a complete contraction
	$$CB(A(G)) \prt A(G) \rightarrow A(G),\; T\otimes f \mapsto T(f),$$
the pointwise multiplication
	\begin{equation}\label{pi-1}
	\pi_1 : A_{cb}(G) \prt A(G) \rightarrow A(G),\; f\otimes g \mapsto fg,
	\end{equation}
is a complete contraction which gives a $A_{cb}(G)$-module structure on $A(G)$.
Moreover $VN(G)$, the dual of $A(G)$, can be equipped with the dual module structure as before;
	\begin{equation}\label{pi-infty}
	\pi_\infty : A_{cb}(G) \prt VN(G) \rightarrow VN(G),\; f \otimes L_{g} \mapsto L_{fg} \ \ (f\in A_{cb}(G), g\in L^1(G)).
	\end{equation}
Clearly $C^*_r(G)$, $UCB(\widehat{G})$, and $C^*_\delta(G)$ are completely bounded $A_{cb}(G)$-submodules of $VN(G)$.

The $A_{cb}(G)$-module structure on $C^*(G)$ can be understood similarly.
We start with a complete contraction
	$$\Phi : CB(B(G)) \prt B(G) \rightarrow B(G),\; T\otimes f \mapsto T(f).$$
Since $A(G)$ is completely contractively complemented in $B(G)$, we have a natural completely isometry $CB(A(G)) \hookrightarrow CB(B(G))$.
If we let
	$$\Psi : A_{cb}(G) \hookrightarrow CB(A(G)) \hookrightarrow CB(B(G)),$$
then we get a complete contraction
	$$\Phi \circ (\Psi \otimes id_{B(G)}) : A_{cb}(G) \prt B(G) \rightarrow B(G),\; f\otimes g \mapsto fg.$$
Indeed, it is easy to check that $\Phi \circ (\Psi \otimes id_{B(G)}) (f\otimes g) = fg$ for any $f\in A(G)$ and $g\in B(G)$.
Since $(C^*(G))^* = B(G)$, we get a complete contraction
	$$A_{cb}(G) \prt C^*(G) \rightarrow C^*(G),\; f\otimes i(g) \mapsto i(fg) \ \ (f\in A_{cb}(G), g\in L^1(G)).$$
It is clear that $C^*_r(G)$, $UCB(\widehat{G})$, and $C^*(G)$ are essential modules.

Recall that a $C^*$-algebra is said to be {\it residually finite-dimensional} if its finite-dimensional,
irreducible $*$-representations separate its points.

	\begin{thm}\label{T:WA-RF-Proj-C(G)}
	Let $G$ be a discrete group.
	If $G$ is amenable or weakly amenable with residually finite-dimensional $C^*(G)$,
	then $C^*_r(G)$, $UCB(\widehat{G})$, and $C^*(G)$ are operator projective in $A_{cb}(G)$-{\bf mod}.
	\end{thm}
\begin{proof}
When $G$ is amenable, we have $A(G) = A_{cb}(G)$, so that the result follows from \cite[Theorem 5.1]{FLS}.

Now suppose that $G$ is a weakly amenable discrete group with residually finite-dimensional $C^*(G)$.
By \cite[Corollary 2.11]{FRS07}, $A_{cb}(G)$ is an operator biprojective, completely contractive Banach algebra with a BAI.
Thus Theorem \ref{thm-essential} leads us to the conclusion.
\end{proof}

\begin{rem}{\rm
Note that $\mathbb{F}_2$ is weakly amenable and $C^*(\mathbb{F}_2)$ is residually finite-dimensional.
}
\end{rem}

	\begin{prop}\label{prop-trivial-mapping-space}
	Let $G$ be a non-discrete group, and $X=C^*_r(G), C^*(G)$ or $C^*_\delta(G)$. Then
		$$_{A_{cb}(G)}B(X, A_{cb}(G)) = 0.$$
	\end{prop}
\begin{proof}
For $X = C^*_r(G)$ or $C^*(G)$ the same proof as \cite[Proposition 5.2]{FLS} can be applied since $A_{cb}(G)\cap L^1(G)$ is dense in $L^1(G)$
and we have a completely contractive embedding
	$$A_{cb}(G) \hookrightarrow L^\infty(G).$$
For $X = C^*_\delta(G)$ we can apply the same proof as \cite[Proposition 7.9]{FLS}.
	
\end{proof}

We also need the following transference result.
	\begin{lem}\label{lem-opensubgp}
	Let $H$ be an open subgroup of a locally compact group $G$.
	If $C^*_r(G)$ (resp. $C^*_\delta(G)$, $C^*(G)$, $UCB(\widehat{G})$, and $VN(G)$) is operator projective in $A_{cb}(G)$-{\bf mod},
	then $C^*_r(H)$ (resp. $C^*_\delta(H)$, $C^*(H)$,  $UCB(\widehat{H})$, and $VN(H)$) is operator projective in $A_{cb}(H)$-{\bf mod}.
	\end{lem}
\begin{proof}
We can apply the proof of \cite[Lemma 5.3]{FLS}. The only ingredients we need more are the following complete contractions.
	$$R_{A_{cb}} : A_{cb}(G) \rightarrow A_{cb}(H), \; g \mapsto g|_H,$$
	$$R_{\text{full}} : C^*(G) \rightarrow C^*(H),\; i(g)\mapsto i(g|_H),$$
and
	$$j_{\text{full}} : C^*(H) \hookrightarrow C^*(G),\; i(f) \mapsto i(\widetilde{f}),$$
where $\widetilde{f}$ is the extension of $f$ to $G$ by assigning 0 outside of $H$.
Indeed, $R_{A_{cb}}$ and $R_{\text{full}}$ are completely contractive since $1_{H}$ is a positive definite function in $M_{cb}A(G)$ with norm 1.
Moreover, $j_{\text{full}}$ is a $*$-homomorphism since it is the extension of an isometric $*$-homomorphism
$L^1(H) \hookrightarrow L^1(G),\; f \mapsto \widetilde{f}$.

\end{proof}

	\begin{thm}\label{T:non-dis. VN(G)-proj}
	Let $G$ be a non-discrete locally compact group. Then we have the following:
		\begin{enumerate}
			\item $C^*_r(G)$ and $C^*(G)$ are not operator projective in $A_{cb}(G)$-{\bf mod}.
			\item When $G$ has an open weakly amenable subgroup, $C^*_\delta(G)$, $UCB(\widehat{G})$ and $VN(G)$ are not operator projective in $A_{cb}(G)$-{\bf mod}.
		\end{enumerate}
	\end{thm}
\begin{proof}
(1) We follow the proof of \cite[Theorem 5.4]{FLS}.
From the structure theory (\cite[Proposition 12.2.2]{Pa01})
we know that there is an almost connected open subgroup $H$ of $G$.
Since $G$ is non-discrete, $H$ is also non-discrete,
and since $H$ is almost connected, $C^*(H)$ is nuclear (\cite{Pat}), which implies $C^*(H)$ has OAP.
Now we suppose $C^*(G)$ is operator projective in $A_{cb}(G)$-{\bf mod}.
Then by Lemma \ref{lem-opensubgp}, $C^*(H)$ is operator projective in $A_{cb}(H)$-{\bf mod}.
By Proposition \ref{prop-proj}, for any non-zero element $x\in C^*(H)$,
we can find a map $T \in {}_{A_{cb}(H)} CB(C^*(H), A_{cb}(H)_+)$ such that $T(x)\neq 0$.
By multiplying an appropriate non-zero function in $A_{cb}(H)$
we can actually find $T' \in {}_{A_{cb}(H)} CB(C^*(H), A_{cb}(H))$ such that $T'(x)\neq 0$.
However, this is impossible by Proposition \ref{prop-trivial-mapping-space}.

For the case of $C^*_r(G)$ we repeat the above argument.
Note that since the nuclearity of a $C^*$-algebra passes to any quotient by a closed ideal \cite[Cor 9.4.4]{BO08},
$C^*_r(H)$ is also nuclear for any almost connected open subgroup $H$ of $G$.

\vspace{0.3cm}
(2) Let $H$ be an open weakly amenable subgroup of $G$. Then $H$ is non-discrete
and $A_{cb}(H)$ has OAP by Corollary \ref{cor-OAP-Acb(G)}. For $C^*_\delta(G)$ we apply the same argument as above.
For $VN(G)$ We fix an open subset $K$ of $H$ with compact closure.
Then it can be shown that $T(L_{1_K}) = 0$ for any $T \in {}_{A_{cb}(H)} B(VN(H), A_{cb}(H))$, as in the proof of \cite[Proposition 5.2]{FLS},
which implies $VN(G)$ is not operator projective in $A_{cb}(G)$-\m\ by Lemma \ref{lem-opensubgp} and Proposition \ref{prop-proj}.
For the case $UCB(\widehat{G})$, we choose $K$ as above, and we choose a nonzero $h\in A_{cb}(H)$ supported in $K$.
Then it can be shown that $T(h\cdot L_{1_K}) = 0$ for any $T \in _{A_{cb}(H)} B(UCB(\widehat{H}), A_{cb}(H))$,
which implies $UCB(\widehat{H})$ is not operator projective in $A_{cb}(H)$-\m. Therefore
$UCB(\widehat{G})$ is not operator projective in $A_{cb}(G)$-\m\ by Lemma \ref{lem-opensubgp}.
\end{proof}

We note that the class of groups satisfying the condition (2) in \ref{T:non-dis. VN(G)-proj}
includes the non-discrete, locally compact groups whose connected component of the identity are
amenable and semisimple connected Lie groups of real rank 1 (see \cite{CDSW} for more details
and examples).

\begin{lem}\label{lem-proj}
	Let $G$ be a weakly amenable discrete group. Suppose that $VN(G)$ is operator projective in $A_{cb}(G)$-{\bf mod}.
	Then there is a bounded projection from $VN(G)$ onto $C^*_r(G)$.
	\end{lem}
\begin{proof}
When $G$ is discrete it is clear that $A_{cb}(G) \cdot VN(G) \subseteq C^*_r(G)$. Since $G$ is weakly amenable we have a BAI in $A_{cb}(G)$.
Then the rest of the proof is the same as \cite[Lemma 3.2]{DP04}.
\end{proof}

If we combine Lemma \ref{lem-proj} with \cite[Theorem 5.7]{FLS}, then we get the following:
	\begin{thm}\label{thm-VN(G)}
	Let $G$ be a discrete group containing an infinite weakly amenable subgroup.
	Then $VN(G)$ is not operator projective in $A_{cb}(G)$-{\bf mod}.
	\end{thm}

\subsection{The modules $L^p(G)$ ($1\le p\le \infty$)}

The pointwise multiplication still gives an $A_{cb}(G)$-module structure on $L^p(G)$ $(1\le p \le \infty)$,
which means that we have the following complete contraction for $1\le p \le \infty$.
	\begin{equation}\label{Lp(G)Acb(G)-mod}
	A_{cb}(G) \prt L^p(G) \rightarrow L^p(G),\; f\otimes g \mapsto fg
	\end{equation}
Indeed, if we combine \eqref{com-Lp} and the following completely contractive inclusion
	$$A_{cb}(G) \hookrightarrow L^\infty(G),$$
we get \eqref{Lp(G)Acb(G)-mod}.

	\begin{thm}
	If $G$ is a discrete, amenable group or a discrete, weakly amenable group with residually finite-dimensional $C^*(G)$,
	then $\ell^p(G)$ $(1\le p< \infty)$ is operator projective in $A_{cb}(G)$-{\bf mod}.
	If $G$ is non-discrete, then $L^p(G)$ $(1\le p< \infty)$ is not operator projective in $A_{cb}(G)$-{\bf mod}.
	\end{thm}

    \begin{proof}
    The proof of the first part is similar to that of Theorem \ref{T:WA-RF-Proj-C(G)}.
    If $G$ is non-discrete, then similar to the argument made in the proof of Proposition
    \ref{prop-trivial-mapping-space}, it follows that $_{A_{cb}(G)}B(L^p(G), A_{cb}(G)) = 0$
    for $1\le p< \infty$. Thus $L^p(G)$ is not operator projective in $A_{cb}(G)$-{\bf mod}
    by Proposition \ref{prop-proj}.
    \end{proof}

For $L^1(G)$ we have a better result.
	\begin{thm}\label{thm-comL1}
	Let $G$ be a locally compact group. If $G$ is discrete,
	then $\ell^1(G)$ is operator projective in $A_{cb}(G)$-{\bf mod}.
	\end{thm}
\begin{proof}
As in the proof of \cite[Theorem 7.2]{FLS} we consider
	$$\rho : \ell^1(G) \rightarrow A_{cb}(G) \prt \ell^1(G),\; \delta_s \mapsto \delta_s \otimes \delta_s.$$
It is straightforward to check that $\rho$ is a left $A_{cb}(G)$-module map, which is a right inverse of the multiplication map.
By the same calculation as in the proof of \cite[Theorem 7.2]{FLS},
we can show that $\rho$ is completely contractive since $\norm{\delta_s}_{A_{cb}(G)} = 1$ for any $s\in G$.
\end{proof}

When $p=\infty$ we have the similar result as in the case of $A(G)$ (\cite[Lemma 7.3]{FLS}) with the same proof .
	\begin{thm}\label{thm-ell-infty}
	Let $G$ be a discrete group. Then $\ell^\infty(G)$ is operator projective in $A_{cb}(G)$-{\bf mod} if and only if $G$ is finite.
	\end{thm}

\subsection{The modules $L^p(VN(G))$ ($1\le p< \infty$)}

As in \cite[section 6]{FLS} and subsection \ref{subsec-LpOSS} we will use complex interpolation. Therefore the case of $p=2$ is important.
So the problem is whether
	\begin{equation*}\label{L2case1}
	\pi_2 : A_{cb}(G)\prt L^2(VN(G)) \rightarrow L^2(VN(G)),\; f \otimes \Lambda(L_{\check{g}}) \mapsto \Lambda(L_{\check{f}\check{g}})
	\end{equation*}
and
	\begin{equation*}\label{L2case2}
	\pi'_2 : A_{cb}(G)\prt L^2(VN(G)) \rightarrow L^2(VN(G)),\; f \otimes \Lambda(L_g) \mapsto \Lambda(L_{fg})
	\end{equation*}
are completely bounded. If we combine the case $p=2$ in \eqref{com-Lp}, a completely contractive inclusion
	$$A_{cb}(G) \hookrightarrow L^\infty(G),$$
and a complete isometry
	$$L^2(VN(G)) \rightarrow L^2(G),\; \Lambda(L_g) \mapsto g,$$
we can conclude that $\pi'_2$ is a complete contraction. The case of $\pi_2$ can be shown similarly.

The proof of the following theorem is similar to that of Theorem \ref{T:WA-RF-Proj-C(G)}.

	\begin{thm}
	Let $G$ be a discrete group.
	If $G$ is amenable or weakly amenable with residually finite-dimensional $C^*(G)$,
	then $L^p(VN(G))$ ($1\le p< \infty$) is operator projective in $A_{cb}(G)$-{\bf mod}.
	\end{thm}

	\begin{prop}\label{prop-trivial-mapping-space2}
	Let $G$ be a non-discrete group. Then for $2\le p<\infty$
		$$_{A_{cb}(G)}B(L^p(VN(G)), A_{cb}(G)) = 0.$$
	\end{prop}
\begin{proof}
The same proof as \cite[Proposition 6.8]{FLS} can be applied since $A_{cb}(G)\cap L^{p'}(G)$ is dense in $L^{p'}(G)$
for $\frac{1}{p} + \frac{1}{p'} = 1$ and we have a completely contractive embedding
	$$A_{cb}(G) \hookrightarrow L^\infty(G).$$
\end{proof}

We also need the following transference result.
	\begin{lem}\label{lem-opensubgp2}
	Let $H$ be an open subgroup of a locally compact group $G$.
	If $L^p(VN(G))$ is operator projective in $A_{cb}(G)$-{\bf mod},
	then $L^p(VN(H))$ is operator projective in $A_{cb}(H)$-{\bf mod}.
	\end{lem}
\begin{proof}
We can apply the proof of \cite[Lemma 6.7]{FLS}. The only ingredient we need more is the following complete contraction.
	$$R_{A_{cb}} : A_{cb}(G) \rightarrow A_{cb}(H), \; g \mapsto g|_H.$$
\end{proof}

	\begin{thm}
	Let $G$ be a non-discrete locally compact group. Then $L^p(VN(G))$ is not operator projective in $A_{cb}(G)$-{\bf mod} for $2\le p <\infty$.
	\end{thm}
\begin{proof}
We follow the proof of \cite[Theorem 5.4]{FLS}.
From the structure theory (\cite[Proposition 12.2.2]{Pa01})
we know that there is an almost connected open subgroup $H$ of $G$.
Since $G$ is non-discrete, $H$ is also non-discrete,
and since $H$ is almost connected, $VN(H)$ is injective (\cite{Pat}), which implies $L^p(VN(H))$ has OAP \cite[Theorem 1.1]{JRX05}.
Now we suppose $L^p(VN(G))$ is operator projective in $A_{cb}(G)$-{\bf mod}.
Then by Lemma \ref{lem-opensubgp2}, $L^p(VN(H))$ is operator projective in $A_{cb}(H)$-{\bf mod}.
By Proposition \ref{prop-proj} for any non-zero element $x\in L^p(VN(H))$
we can find a map $T \in {}_{A_{cb}(H)} CB(L^p(VN(H)), A_{cb}(H))$ such that $T(x)\neq 0$.
However, this is impossible by Proposition \ref{prop-trivial-mapping-space2}.
\end{proof}

\section{Summary of the results}

We will end the paper with tables that contain a summary of our results.
The first column in the table denotes the module under consideration.
The second column identifies those classes of groups for which we know 
definitively that the module under consideration is operator projective as a left $A$-module for the corresponding algebra $A$.
The third column identifies those classes of groups for which we know 
definitively that the module under consideration is not operator projective as a left $A$-module for the corresponding algebra $A$.
Finally, the question mark in the second table implies we could not draw any conclusion for the question.

\begin{table}
\caption{Summary of the results}
\begin{tabular}{|p{4cm}|p{2cm}|p{4cm}|}

\hline $A$-{\bf mod} for $A=S^1A(G)$ or $S_0(G)$ & op. proj. & not op. proj. \\

\hline
\hline $A(G)$ & (1) & $G$ is discrete and infinite or $G = SL(3,\mathbb{R})$ \\

\hline $C^*_r(G)$, $UCB(\widehat{G})$  & (1) & (2) \\

\hline $L^p(VN(G))$ $(1<p<\infty)$ & (1) & (2) \;\;\; when \;\;\;\;$2\le p <\infty$  \\

\hline $VN(G)$ & (1) & (2)  \\

\hline $C^*(G)$, $C^*_\delta (G)$ & (1) & (2)  \\

\hline $L^1(G)$ & $G$ is discrete & $G$ is non-discrete \\

\hline $L^p(G)$ $\;\;\;(1<p<\infty)$ & (1) & (2)  \\

\hline $L^\infty(G)$ & (1) & $G$ is discrete and infinite  \\

\hline $A(G)^{**}$, $UCB(\widehat{G})^*$ & (1) & (2) \\

\hline \multicolumn{3}{l}{(1) : $G$ is finite.}\\

\hline \multicolumn{3}{l}{(2) : $G$ is infinite.}\\

\hline

\end{tabular}

\end{table}

\begin{table}
\caption{Summary of the results}
\begin{tabular}{|p{3.9cm}|p{1.9cm}|p{3.9cm}|}

\hline $A_{cb}(G)$-{\bf mod} & op. proj. & not op. proj.  \\

\hline
\hline $A(G)$ & [IN] & ? \\

\hline $C^*_r(G)$ & (1) & (2) \\

\hline $L^p(VN(G))$ $(1<p<\infty)$ & (1) & (2) \;\;\; when \;\;\;\;$2\le p <\infty$ \\

\hline $VN(G)$ & $G$ is finite & (3) or (4) \\

\hline $C^*(G)$ & (1) & (2) \\

\hline $L^1(G)$ & $G$ is discrete & (2) \\

\hline $L^p(G)$ $\;\;\;(1<p<\infty)$ & (1) & (2) \\

\hline $L^\infty(G)$ & $G$ is finite & $G$ is infinite and discrete \\

\hline $C^*_\delta (G)$, $UCB(\widehat{G})$ & (1) & (3) \\

\hline $A(G)^{**}$, $UCB(\widehat{G})^*$ & $G$ is finite & $G$ is infinite and amenable \\

\hline \multicolumn{3}{l}{(1) :$\begin{array}{l} \text{$G$ is discrete, weakly amenable with residually finite-dimensional $C^*(G)$}\\
\text{or $G$ is discrete, amenable.}\end{array}$}\\

\hline \multicolumn{3}{l}{(2) : $G$ is non-discrete.}\\

\hline \multicolumn{3}{l}{(3) : $G$ is non-discrete with a weakly amenable open subgroup.}\\

\hline \multicolumn{3}{l}{(4) : $G$ is discrete group containing an infinite weakly amenable subgroup.}\\

\hline
\end{tabular}

\end{table}

\newpage

\bibliographystyle{amsplain}
%\bibliography{hidden}
\providecommand{\bysame}{\leavevmode\hbox to3em{\hrulefill}\thinspace}

\end{document}